\documentclass[journal,twoside,web]{ieeecolor}
\usepackage{generic}
\usepackage{cite}
\usepackage{amsmath,amssymb,amsfonts}
\usepackage{algorithmic}
\usepackage{graphicx}
\usepackage{textcomp}
\newtheorem{lemma}{Lemma}
\newtheorem{definition}{Definition}
\newtheorem{theorem}{Theorem}
\newtheorem{remark}{Remark}

\newtheorem{Assumption}{Assumption}
\def\BibTeX{{\rm B\kern-.05em{\sc i\kern-.025em b}\kern-.08em
    T\kern-.1667em\lower.7ex\hbox{E}\kern-.125emX}}

\begin{document}
\title{The topology in the game controllability of multiagent systems}
\author{Junhao Guo, Zhijian Ji, and Yungang Liu
\thanks{This work was supported by the National Natural Science Foundation of China(Grant. 61873136, 62033007), Taishan Scholars Climbing Program of Shandong Province of China and Taishan Scholars Project of Shandong Province of China(Grant. ts20190930).}
\thanks{Junhao Guo, Zhijian Ji  and Chong Lin are with Shandong Key Laboratory of Industrial Control Technology, Institute of Complexity Science, College of Automation, Qingdao University, Qingdao, Shandong, 266071, China.(941708869@qq.com;jizhijian@pku.org.cn)}
\thanks{Yungang Liu is with School of Control Science and Engineering, Shandong University, Jinan, 250061, China.(lygfr@sdu.edu.cn)}
}
\maketitle
\begin{abstract}
  In this paper,  the graph-based condition for the controllability of game-based  control system is presented when the control of regulator is not zero. A  control framework which can describe realism well--expressed as the game-based control system (GBCS),  was obtained  in 2019, which, unfortunately, is not graph theoretically verifiable, and the   regulator's control input is assumed to be zero.
 However, based on a new established notion, strategy matrix, we propose a graph theory condition to judge the controllability of GBCS, instead of using algebraic conditions for complex mathematical calculations. More specifically,
to tackle these issues, one needs to  study the expression of Nash equilibrium actions when regulator's control is not zero first.
Based on this expression, the general formula of game controllability matrix is obtained, which provides theoretical support for studying the essential influence of topology on game-based control system.
The general formula is always affected by the specific matrix--strategy matrix, composed of Nash equilibrium  actions, and the matrix can not only be obtained by matrix calculation, but also can be directly written through the topology, which is  the specific influence of the topology   on the GBCS.
Finally, we obtain the result of judging the controllability of the system directly according to the topological structure, and put forward the conjecture that there is no limitation of equivalent partition in GBCS.
Arguably, this is a surprising conjecture on the equivalent partition  of graphs, because only the limitation of equivalent partition in five-node graphs has been solved so far.
\end{abstract}

\begin{IEEEkeywords}
 Controllability, game-based control system, topological structure, equivalent partition
\end{IEEEkeywords}

\section{Introduction}
\label{sec:introduction}
\IEEEPARstart{G}{ame} is ubiquitous in nature. In recent years, game theory has gradually become a powerful tool to study control systems, especially for distributed control systems\cite{jason,mar}. This leds to the introduction and research of game-based control systems (GBCS)\cite{zhangrenren}. This system has a hierarchical decision structure: one regulator and multiple agents.
Although the study of GBCS has some practical significance, when the game is defined on a large-scale system, each player's strategy often depends on the structure of the underlying network. Therefore, it is of great significance to establish a game control model related to topology structure.
Going through the literature, one can find many algebraic conditions  for judging controllability of the GBCS.
 However, when faced with massive game-based control system, it is very complicated to judge its controllability by mathematical calculation. This  makes us establish the graph-based criterion of the game-based control system.
It is natural then to ask  what effect does topology have on the game-based control  system? It is found that it has an important influence on the formation of Nash equilibrium actions.

\subsection{Literature review}
The game-based control system originates from the regulator problem of linear time-invariant systems proposed by Jacob et al. That is, one target may be to find a control function $u(t)$ which drives the state of the system  to a small neighborhood of zero at time $T$ \cite{eng}.
In contrast, if the state of the system represents a set of financial related economic variables, then, $u$ represents the investment behavior that can increase these variables with the goal of controlling these variables to reach the desired level as quickly as possible.
Zhang  et al. \cite{zhangrenren} studied a linear   system with a hierarchical decision structure consisting of a regulator and multiple agents.
  The regulator first makes a decision, and then each agent optimizes its own payoff function to reach the possible Nash equilibrium of the non-cooperative dynamic game.
 At this time, the underlying agents in the system have information transfer with the regulator.
This can be reflected in the topology structure, that is,  regulator changes the macrostate depending on Nash equilibrium actions which is formed based on the characteristics of topology structure. If the positions of the underlying agents change, the Nash equilibrium actions will also change, which will affect the controllability of the system.

At present, the influence of network topology structure on the  controllability of multiagent systems (MASs ) is mainly studied by graph theory, as shown in [5-10]. It can be seen that the controllability of the MASs is closely related to the underlying graph topology, and graph division\cite{car,god} describes the controllability of the MASs from the perspective of topology.
Cardoso et al. proposed a sufficient and necessary condition for almost equivalent partition, and clarified the relationship between Laplacian matrix and general Laplacian matrix \cite{car}. Qu et al. investigated the controllability of multi-agent systems under equivalent partition  \cite{quji}.
After devoting much thought to a graph-theoretical problem, one gradually realizes that it is a new trend to use graph theory to study game behavior in multi-agent system[14-15]. For instance, the pursuit-evasion game problem \cite{xiang1}, which   has been extensively used in aircraft control and missile guidance, intelligent transportation system and collision avoidance design of wireless sensor network in military implementation. Using graph theory to study the interaction between agents with limited perceptive ability, so as to obtain the distributed control strategy of each agent, has gradually become the key to study this problem. The establishment of Bayesian graph game model is also inseparable from the principle of local interaction and game behavior fusion. In literature \cite{xiang11}, Bayesian games (games with incomplete information) describe an agent participating in an unspecified game where the true intentions of other players may be unknown and each player must adjust its goals accordingly. Lopez et al.  \cite{xiang14} proposed two belief updating methods that do not require graphical topology knowledge. The first method is the application of Bayesian rules, and the second method is treated as a modification of non-Bayesian updates. Mohammad et al. \cite{xiang15} used structured systems theory and other graph theory concepts to analyze games, improved the detectability of network physical attacks, and studied the optimal configuration of sensors in networked control systems. Although the above researches have made great progress, the topological characterization of game systems based on Nash Equilibrium action remains largely unknown.

\subsection{Contributions}
We turn to the more concrete task of discussing the contents of this article:
Compared with the traditional control theory framework, the game-based control system will consider the strategic behavior of each agent, which will avoid the unrealistic phenomenon of system dynamics caused by ignoring the agent's own behavior.
However, the current research on this system is only embodied in the algebraic level.
Whether the neighbor relationship of the agents playing with each other has an impact on this system?
Or, stated differently, if the position relationship between agents is regarded as a network topology, what essential influence does this network topology have on the internal mechanism of game-based control system? The study of this problem has very obvious research significance, that is, when studying the controllability of game-based control system, we do not need to carry out complex mathematical calculations, just need to judge according to the graph theory conditions.

To shed light into this problem, this paper first analyzes the framework of game-based control system.
Different from the previous literature, which assumes that regulator's control strategy is zero, we obtain the Nash equilibrium action expression of the system under regulator's control.
Based on this expression, we obtain the general formula of  controllability matrix.
We find that the strategy matrix, which plays a decisive role in the general formula, can not only be obtained by algebraic calculation, but also can be written directly according to the topology structure, which is the fundamental influence of the topology structure on the game-based control system.

Under the premise that game-based control system has become a big research trend, the strategy matrix may become an important research tool to replace the Laplacian matrix.
Finally, a graph-theoretic condition for the controllability of game-based control systems is obtained. Compared with the general system that does not consider game factors, it also requires special consideration of the relationship between the row vectors of the matrix that corresponds to the multiplication of  the symmetric solution of the Riccati equation and  the transition matrix.
\subsection{Notation}
In this paper, different agents are regarded as nodes and their communication relations as edges. The topological relation of the system can be expressed by graph $G$, which consists of node set $V(G)$ and edge set $E(G)\subset V(G)\times V(G)$. Assuming that there are one regulator and $H$ bottom agents in the system, then $V(G)=\{r_{1},l_{1},...,l_{H}\}, E(G)=\{(i,j)|i,j\in V(G)\}$, and the neighbor of node $l_{i}$ is defined as $N(l_{i})=\{(l_{i},l_{j})|l_{i},l_{j}\in E(G)\}$. Graph $G$ is connected if there is a path between any two different nodes $i$ and $j$ in graph $G$. The adjacency matrix of graph $G$ is defined as
$$
 A(G)=[a_{ij}]_{n\times n}=\left \{
\begin{aligned}
1~~(i,j)\in E(G)  \\
0~~(i,j)\notin E(G)\\
\end{aligned}
\right..
$$
The Laplacian matrix of the graph $G$ is defined as $L(G)=D(G)-A(G)$, where $D(G)=diag([d_{i}]^{n}_{i=1})$, $d_{i}$ indicates the number of neighbors of node $i $.
$I_{n}$ is the $n\times n$-dimensional identity matrix, and $0_{n\times m}$  denotes the null matrix of $n\times m$-dimension (or $0_{n}$, if $n= m$ ). $1_{n}=[1,1,...,1]^{T}$ represents an $n$-dimensional column vector of $n$ ones.

\subsection{Organization}
The remainder of this paper is organized as follows. In Section II,  we introduce the game-based control system and analyze the representation of topology structure of the system;
 Section III gives the expression of Nash equilibrium actions under the premise that the regulator's control  input is not zero.
The graph theory condition of game-based control system is given in Section IV. Section V concludes this paper.

\section{Game-based control system}
\label{sec:guidelines}

\subsection{Model analysis }
The dynamics for a continuous-time linear time-invariant system is
\begin{equation}\label{gs1}
  \dot{x}(t)=Ax(t)+Bu(t),
\end{equation}
where $x(t)=[x_{1}(t),\dots, x_{n}(t)]^{T}$ is the stacked system states, $u(t)\in R^{m}$ denotes the control input.  When studying the influence of topology structure on the controllability of the system, all agents in the system are divided into the set of leaders and the set of followers. At this time, the form of control input is based on the topological structure relationship formed between leaders and followers, and the system matrix $A$ is represented by the Laplacian matrix $L$.

Based on the idea of a two-player zero-sum game, the following linear differential equation describes the dynamics of the game with $H$ players\cite{engwerda}:
\begin{equation}\label{gs2}
  \dot{x}(t)=Ax(t)+B_{1}u_{1}(t)+...+B_{H}u_{H}(t).
\end{equation}
where $x(t)$ still represents the state of the system, and $ u_{i},i=1,2, \dots, H,$  represent the  $m_{i}$ -dimensional vector that can be manipulated by player $i$. Each player aims to minimize their quadratic cost function
\begin{equation}\label{gs3}
\begin{aligned}
  J_{i}(u_{1},...,u_{H})=&\int_{0}^{T}\{x^{T}(t)Q_{i}x(t)+\sum_{j=1}^{H}u^{T}_{j}(t)R_{ij}u_{j}(t)\}dt
\\&+x^{T}(T)Q_{iT}x(T),i=1,2,...H.
\end{aligned}
\end{equation}
Based on the form of systems (\ref{gs1}) and (\ref{gs2}), GBCS  with one regulator and $H$ agents can be expressed as:
\begin{equation}\label{gs4}
\begin{footnotesize}
\begin{aligned}
  \dot{x}_{r}(t)=&Ax_{r}(t)+\sum_{i=1}^{H}A_{i}x_{i}(t)+\sum_{i=1}^{H}D_{i}u_{i}(t)
+Bu(t), x(0)=x_{0}.
\end{aligned}
\end{footnotesize}
\end{equation}
\begin{equation}\label{gs5}
\begin{footnotesize}
\begin{aligned}
  &\dot{x}_{i}(t)=E_{i}x_{r}(t)+\sum_{j=1}^{H}F_{ij}x_{i}(t)+\sum_{j=1}^{H}B_{ij}u_{j}(t)
+B_{i}u(t),\\& x_{i}(0)=x_{i,0}, i=1,2,...,H.
\end{aligned}
\end{footnotesize}
\end{equation}
The payoff function to be minimized by $u_{i}(\cdot)$ of any agent $i(i =1, 2, . . . , H)$ is
\begin{equation}\label{gs6}
\begin{aligned}
  J_{i}(u_{1},...,u_{H})=&\frac{1}{2}\int_{0}^{T}\{X^{T}(t)Q_{i}X(t)+u^{T}_{i}(t)R_{i}u_{i}(t)\}dt
\\&+\frac{1}{2}x^{F}(T)^{T}Q_{iT}x^{F}(T),i=1,2,...H.
\end{aligned}
\end{equation}
In system (\ref{gs1}), the influence of topology structure on the controllability of the system is reflected in the system matrix $L$, for the game-based control system  (\ref{gs4}) (\ref{gs5}) (\ref{gs6}), the influence of topology structure on the system is not only reflected in the system matrix, but also reflected in each action $u_{i}$. Therefore, this paper studies the influence of topology structure on GBCS under the following conditions.

\subsection{Model description}
Consider a system with one regulator and $H$ agents

$
\begin{aligned}
\dot{x}_{r}(t)=ax_{r}(t)+\sum_{i=1}^{H}a_{i}x_{i}(t)+\sum_{i=1}^{H}d_{i}u_{i}(t)+bu(t),   
\end{aligned}
$
\begin{equation}\label{gs7}
\begin{aligned}
    &\dot{x}_{i}(t)=e_{i}x_{r}(t)+\sum_{i=1}^{H}f_{ij}x_{i}(t)+\sum_{i=1}^{H}b_{ij}u_{j}(t)
+b_{i}u(t),   
   \\& x(0)=x_{0},x_{i}(0)=x_{i,0}.
\end{aligned}
\end{equation}
Each player has a quadratic cost function:
\begin{equation}\label{gs8}
\begin{aligned}
  J_{i}(u_{1},...,u_{H})=&\frac{1}{2}\int_{0}^{T}\{X^{T}(t)Q_{i}X(t)+ \sum_{j=1}^{H}u^{T}_{j}(t)R_{ij}u_{j}(t)\}dt\\&+\frac{1}{2}x^{F}(T)^{T}Q_{iT}x^{F}(T),i=1,2,...H.
\end{aligned}
\end{equation}
where $x^{F}=\left(
           \begin{array}{c}
             x_{1} \\
               x_{2} \\
             \vdots \\
             x_{H} \\
           \end{array}
         \right), X(t)=\left(
                         \begin{array}{c}
                           x_{r}(t) \\
                            x^{F}(t) \\
                         \end{array}
                       \right)$. $Q_{i}, Q_{iT}$ are symmetric matrices, $R_{i}>0,$
Let
\begin{equation*}\label{gs9}
 \tilde{A}=\left(
               \begin{array}{ccccc}
                 a & a_{1} & a_{2} &\dots & a_{H} \\
                  e_{1} & f_{11} & f_{12} & \dots &  f_{1H} \\
                 e_{2} &  f_{21} &  f_{22} & \dots & f_{2H} \\
                 \vdots & \vdots & \vdots &  & \vdots \\
                  e_{H} &  f_{H1} &  f_{H2}&\dots &  f_{H} \\
               \end{array}
             \right),
\end{equation*}
\begin{equation}
 \tilde{B}_{i}=\left(
                                       \begin{array}{c}
                                         d_{i} \\
                                         b_{1i} \\
                                          b_{2i} \\
                                         \vdots \\
                                          b_{Hi} \\
                                       \end{array}
                                     \right),  \tilde{B}=\left(
                                                           \begin{array}{c}
                                                             b \\
                                                             b_{1} \\
                                                             b_{2} \\
                                                             \vdots \\
                                                             b_{H} \\
                                                           \end{array}
                                                         \right), \bar{B}=\left(
                               \begin{array}{c}
                                 \tilde{B} \\
                                 0 \\
                                 0 \\
                                 \vdots \\
                                 0 \\
                               \end{array}
                             \right).
\end{equation}

$ \bar{A}=$
$$ \begin{footnotesize}
\left(
            \begin{array}{ccccc}
              \tilde{A} & \tilde{B}_{1}R_{1}(\tilde{B}_{1}+C)^{T} &  \tilde{B}_{1}R_{1}(\tilde{B}_{1}+C)^{T} & \dots &  \tilde{B}_{1}R_{1}(\tilde{B}_{1}+C)^{T} \\
              Q_{1} & - \tilde{A}^{T} &0 &  \dots  & 0 \\
              Q_{2} & 0 &  - \tilde{A}^{T} & \dots & 0 \\
              \vdots &   \vdots &  \vdots   &   &   \vdots \\
              Q_{L} & 0 & 0 &  \dots  &  - \tilde{A}^{T} \\
            \end{array}
          \right).
\end{footnotesize}$$
In order to more clearly show the essential influence of topology structure on game-based control system, we only consider the case where $x(t), x_{i}(t), u_{i}(t),u(t)$  is one-dimensional, and the corresponding $a, a_{i}, d_{i},b, e_{i}, f_{ij}, b_{ij}$ and $ b_{i}$  are all constant.  Therefore, $\bar{A}\in R^{(1+H)^{2}\times (1+H)^{2}}, \bar{B}\in R^{(1+H)^{2}}$.
As the elements in $\bar{A}$ are mostly studied in the form of block matrix in the following, each block matrix $\tilde{A}, Q_{i} $ and $\tilde{B}_{i}R_{i}^{-1}(\tilde{B}_{i}+C)^{T}$ in $\bar{A}$  is regarded as a whole, and $\bar{A}_{ij}$ represents the block matrix in row $i$ and column $j$.

The idea expressed in this model is: The regulator makes the decision first, and the underlying agents reach Nash equilibrium after receiving the decision of the regulator, so the state of the system is essentially caused by regulator's control. Therefore, the game strategy formed by $H$ agents is not in parallel with the external control input, but according to regulator's input, Nash equilibrium action is formed.
When the regulator's strategy is not zero, the Nash equilibrium action of the system is unknown, so we first solve for the Nash equilibrium action.

\section{Nash equilibrium action of the GBCS }
By comparing the above different models, it can be found that system  (\ref{gs4})(\ref{gs5}) adds the term $[B_{1}u_{1}(t)+...+B_{H}u_{H}(t)]$ in system (\ref{gs2}) to the original control system (\ref{gs1}), for better differentiation, (\ref{gs7}) can be written as
\begin{equation}\label{gs11}
  \dot{x}(t)=Ax(t)+\sum_{i=1}^{H}B_{i}u_{i}(t)+Cz(t).
\end{equation}
So the difference between (\ref{gs11}) and (\ref{gs1})  is the action $u_{i}$.  Here, $u_{i}$ is the action  made by agent $i$ after the non-cooperative differential game. If we consider the topological relationship between agents, for example, agents $i$ and $j$ are not neighbors, it can be seen from  (\ref{gs7}) that for $\dot{x}_{i}(t)$, the coefficient in front of $u_{j}$  can be zero. Similarly, for $\dot{x}_{j}(t)$, the coefficient in front of $u_{i}$ is zero. Therefore, the coefficient $B_{i}$ in front of $u_{i}$  is an important index reflecting the influence of the topology on the game. The Nash equilibrium action $u^{*}_{i}$ reached by players after the game is related to the coefficient $B_{i}$, that is, the Nash equilibrium action obtained in system (\ref{gs2}) (\ref{gs3})  is $u^{*}_{i}=-R_{i}^{-1}B_{i}^{T}(t)\varphi_{i}(t)$. However, the condition for this Nash equilibrium action $u^{*}_{i}$ to hold under system  (\ref{gs7}) is $z(t)=0$ (\cite{zhangrenren}Assumption 1). Therefore, the controllability conditions of GBCS obtained at present are based on the assumption that $z(t)=0$, however, this assumption indicates that regulator does not play the role of control. Therefore, if the complete control system is considered, the assumption that $z(t)=0$ cannot be based. When $z(t)$ is not equal to 0, the Nash equilibrium action  $u^{*}_{i}=-R_{i}^{-1}B_{i}^{T}(t)\varphi_{i}(t)$ obtained in \cite{zhangrenren} is not necessarily valid, because the conclusion is obtained based on system  (\ref{gs2}) (\ref{gs3}). Therefore, this section studies the Nash equilibrium action based on  (\ref{gs11}) .

\lemma Consider the optimal problem of game-based control systems :

\begin{equation}\label{gs12}
  J(x_{0}, u)=\int_{0}^{T}g(t,x(t),u(t))dt+h(x(T)),
\end{equation}
\begin{equation}\label{gs13}
   \dot{x}(t)=f(t,x(t),u(t),z(t)), x(0)=x_{0}.
\end{equation}
Let $H(t,x(t),u(t), \lambda(t),z(t)):=g(t,x(t),u(t))+\lambda(t)f(t,x(t),u(t),z(t))$,  suppose $g(t,x(t),u(t))$ and $f(t,x(t),u(t),z(t))$  are continuous;
 For $f, g$, the partial derivative with respect to $x, u$, and $z$ exists and is unique;
 $h(x)$ is a continuously differentiable function.
 If $u^{*}(t)$ is the action that minimizes  (\ref{gs12}) and $x^{*}(t), \lambda^{*}(t)$ and $ z^{*}(t)$ are the corresponding state, parameter and control input respectively, then
\begin{equation}\label{gs14}
\scriptsize
\begin{aligned}
 & \dot{x}^{*}(t)=f(t,x^{*}(t), u^{*}(t),   z^{*}(t))=(=\frac{\partial H(t,x^{*}(t),u^{*}(t), \lambda^{*}(t),z^{*}(t))}{\partial \lambda}), 
\\& x^{*}(0)=x_{0}.
\end{aligned}
\end{equation}
\begin{equation}\label{gs15}
\scriptsize
  \lambda^{*}(t)=-\frac{\partial H(t,x^{*}(t),u^{*}(t), \lambda^{*}(t),z^{*}(t))}{\partial x};
 \lambda^{*}(T)=\frac{\partial h(x^{*}(T))}{\partial x}.
\end{equation}

For all $t\in [0,T]$,
\begin{equation}\label{gs16}
\scriptsize
  \frac{H(t,x^{*}(t),u^{*}(t), \lambda^{*}(t),z^{*}(t))}{\partial u}+\frac{H(t,x^{*}(t),u^{*}(t), \lambda^{*}(t),z^{*}(t))}{\partial z}=0.
\end{equation}

\begin{proof}
According to  (\ref{gs13}), a Lagrange multiplier $\lambda(t)$ is chosen arbitrarily, then
\begin{equation}\label{gs17}
\scriptsize
  \int_{0}^{T} \lambda(t)[f(t,x(t),u(t),z(t))-\dot{x}(t)]dt=0.
\end{equation}
Append  (\ref{gs17}) to  (\ref{gs12})
\begin{equation}\label{gs18}
\scriptsize
   \bar{J}(x_{0}, u)=\int_{0}^{T}\{g(t,x(t),u(t)+\lambda(t)[f(t,x(t),u(t),z(t))-\dot{x}(t)]\}dt+h(x(T)).
\end{equation}
Define the Halmitonian function
\begin{equation}\label{gs19}
\scriptsize
  H(t,x(t),u(t), \lambda(t),z(t)):=g(t,x(t),u(t))+\lambda(t)f(t,x(t),u(t),z(t)),
\end{equation}
then
\begin{equation}\label{gs20}
\scriptsize
  \bar{J}(x_{0}, u)=\int_{0}^{T}\{H(t,x(t),u(t), \lambda(t),z(t))-\lambda(t)\dot{x}(t)]\}dt+h(x(T)).
\end{equation}
Integration by parts shows that
\begin{equation}\label{gs21}
  -\int_{0}^{T}\lambda(t)\dot{x}(t)dt=-\lambda(T)x(T)+\lambda(0)x_{0}+\int_{0}^{T}\dot{\lambda}(t)x(t)dt.
\end{equation}
Hence,
\begin{equation}\label{gs22}
\begin{aligned}
  \bar{J}(x_{0}, u)=&\int_{0}^{T}\{H(t,x(t),u(t), \lambda(t),z(t))-\lambda(t)\dot{x}(t)]\}dt\\&+h(x(T)).
\end{aligned}
\end{equation}
According to   (\ref{gs22}), the last three terms in $\bar{J}(x_{0}, u)$ are only related to the initial time, and are independent of $t$. Therefore, no matter how the path of $\lambda(t)$ is chosen, it will have no influence on the value of $\bar{J}(x_{0}, u)$. That is, $\dot{x}(t)=\frac{\partial H}{\partial \lambda}, t\in[0,T]$  forms the necessary condition for $\bar{J}$ to take an extreme value.
Now assume that $u^{*}$ is the optimal strategy to minimize $\bar{J}(x_{0}, u)$, and $x^{*}(t)$ is the corresponding optimal state trajectory. If $u^{*}$ is slightly disturbed, then $u(t)=u^{*}(t)+\epsilon p(t)$, if $\epsilon$  is small enough,   (\ref{gs22}) becomes
\begin{equation}\label{gs23}
\begin{aligned}
  \bar{J}(\epsilon)=&\int_{0}^{T}[H(t,x(t,\epsilon,p), u^{*}(t)+\epsilon p, \lambda, z(t,\epsilon,p))]dt
\\&+h(x(T,\epsilon,p)-\lambda(T)x(T,\epsilon,p)+\lambda (0)x_{0}.
\end{aligned}
\end{equation}
Since $ \bar{J}(\epsilon)$ minimizes at $\epsilon=0$  and $f,g,h,z$ are differentiable, $\bar{J}$  is differentiable with respect to $\epsilon$. Then, $\frac{dJ}{d\epsilon}=0 $ at $\epsilon=0,$ where
\begin{equation}\label{gs24}
\begin{aligned}
  \frac{d\bar{J}}{d\epsilon}=&\int_{0}^{T}\{\frac{\partial H}{\partial x}\frac{\partial x}{\partial \epsilon}+\frac{\partial H}{\partial u}p(t)+\frac{\partial H}{\partial z}\frac{\partial z}{\partial \epsilon}+\dot{\lambda}(t)\frac{\partial x}{\partial \epsilon}\}dt
\\&+\frac{\partial h(x(T,\epsilon,p))}{\partial x}\frac{dx}{d\epsilon}-\lambda(T)\frac{dx}{d\epsilon}.
\end{aligned}
\end{equation}
Since $\lambda(t)$ is arbitrary, the equality $\frac{d\bar{J}}{d\epsilon}=0 $ at $\epsilon=0$ is still true if we choose the solution of the boundary value problem below as the value of $\lambda(t)$:
$\frac{\partial H(t,x^{*},u^{*},\lambda,z^{*})}{\partial x}+\dot{\lambda}=0, $ with $\frac{\partial h(x^{*}(T))}{\partial x}-\lambda(T)=0.$
Since the partial derivatives of $f,g,z$ with respect to $x$ are continuous, then $H(t,x^{*},u^{*},\lambda,z^{*})=g(t,x^{*},u^{*})+\lambda f(t,x^{*},u^{*},z^{*})$, $\frac{\partial H}{\partial x}=\frac{\partial g}{\partial x}+\lambda\frac{\partial f}{\partial x}$.
 That is $\dot{\lambda}=-\frac{\partial f}{\partial x}\lambda+\frac{\partial g}{\partial x}.$
 So there is a unique solution to the boundary value problem.
 Let $\frac{\partial \bar{J}}{\partial \epsilon}=0$ , $x=x^{*}, u=u^{*}$, then
\begin{equation}\label{gs25}
\begin{aligned}
  \frac{\partial \bar{J}}{\partial \epsilon}=&\int_{0}^{T}\{[\frac{\partial H}{\partial x}+\dot{\lambda}(t)]\frac{dx}{d\epsilon}+\frac{\partial H}{\partial u}p(t)+\frac{\partial H}{\partial z}\frac{dz}{d\epsilon}\}dt
\\&+[\frac{\partial h(x^{*})}{\partial x}-\lambda(T)]\frac{dx(T,\epsilon,p)}{d\epsilon}.
\end{aligned}
\end{equation}
According to the boundary value conditions,
\begin{equation}\label{gs26}
  \int_{0}^{T}\{\frac{\partial H(t,x^{*},u^{*},\lambda,z^{*})}{\partial u}p(t)+\frac{\partial H(t,x^{*},u^{*},\lambda,z^{*})}{\partial z}\frac{dz}{d\epsilon}\}dt=0.
\end{equation}
Since the above equation is true for any continuous function $p(t)$ on $[0,T]$, if $p(t)=\frac{dz}{d\epsilon}$, then
\begin{equation}\label{gs27}
  \frac{\partial H(t,x^{*},u^{*},\lambda,z^{*})}{\partial u}+\frac{\partial H(t,x^{*},u^{*},\lambda,z^{*})}{\partial z}=0.
\end{equation}
\end{proof}

\Assumption{For any initial state $x_{0}, x_{i,0}$, open-loop Nash equilibria exists in (\ref{gs7})(\ref{gs8}), and the following Riccati differential equations have symmetric solutions $K_{i}, i=1,2,\dots,H$
\begin{equation}\label{24}
\begin{aligned}
  &\dot{K}_{i}(t)=-\tilde{A}^{T}K_{i}(t)-K_{i}(t)\tilde{A}+K_{i}(t)\tilde{S}_{i}K_{i}(t)-Q_{i}, \\&K_{i}(T)=\tilde{Q}_{iT}, i=1,2, \tilde{S}_{i}:=\tilde{B}_{i}R_{ii}^{-1}(\tilde{B}_{i}+C)^{T}
\end{aligned}
\end{equation}}

\theorem
Consider a linear system with two players
\begin{equation}\label{gs29}
\begin{aligned}
  J_{i}(u_{1},u_{2}):=&\int_{0}^{T}\{x^{T}(t)Q_{i}x(t)+u_{1}^{T}R_{i1}u_{1}+u_{2}^{T}R_{i2}u_{2}\}dt\\&+x^{T}(T)Q_{i1}x(T),
\end{aligned}
\end{equation}
\begin{equation}\label{gs30}
\dot{x}=Ax+B_{1}u_{1}+B_{2}u_{2}+Cz,
\end{equation}
\begin{equation}\label{gs31}
  M=\left(
      \begin{array}{ccc}
        A & -S_{1} & -S_{1} \\
        -Q_{1} & -A^{T} & 0 \\
        -Q_{2} & 0 & -A^{T} \\
      \end{array}
    \right),  {S}_{i}:= {B}_{i}R_{ii}^{-1}( {B}_{i}+C)^{T}.
\end{equation}
Then the necessary and sufficient condition for the existence of an open-loop Nash equilibrium action for each initial state is that the matrix $H(T)$ is invertible,
\begin{equation}\label{gs32}
  H(T)=\left(
         \begin{array}{cccc}
           I & 0 & 0& 0\\
         \end{array}
       \right)e^{-MT}\left(
                 \begin{array}{c}
                   I \\
                   Q_{1T} \\
                   Q_{2T} \\
                   0 \\
                 \end{array}
               \right).
\end{equation}
Moreover,
if for every $x_{0}$,  there exists an open-loop Nash equilibrium action, then the action is unique, and
$u_{i}=-R_{ii}^{-1}(B_{i}+C)^{T}\varphi_{i}(t)$.

\begin{proof}
(``$\Rightarrow$" part)
If $(u_{1}^{*}, u_{2}^{*})$ is a Nash equilibrium action, then according to Lemma 1, the Hamiltonian function of player $i$ for action $u_{i}$ is
$$
\begin{aligned}
H_{i}(t,x,u_{1},u_{2},\psi_{i},z)=&(x^{T}(t)Q_{i}x(t)+u_{1}^{T}R_{i1}u_{1}+u_{2}^{T}R_{i2}u_{2})
\\&+\psi^{T}_{i}(Ax+B_{1}u_{1}+B_{2}u_{2}+Cz).
\end{aligned}$$
When $J_{i}$ is minimized, $\frac{\partial H_{i}(t,x^{*},u^{*},\lambda^{*},z^{*})}{\partial u_{i}}+\frac{\partial H_{i}(t,x^{*},u^{*},\lambda^{*},z^{*})}{\partial z}=0.$
That is
$$u_{1}^{*}(t)=-R_{11}^{-1}(B_{1}+C)^{T}\psi_{1}(t), u_{2}^{*}(t)=-R_{22}^{-1}(B_{2}+C)^{T}\psi_{2}(t),$$
where $\psi_{i}, i=1,2$ satisfy the following conditions
$$\dot{\psi}_{1}(t)=-Q_{1}x(t)-A^{T}\psi_{1}(t), with ~ \psi_{1}(t)=Q_{1T}x(T),$$
$$\dot{\psi}_{2}(t)=-Q_{2}x(t)-A^{T}\psi_{1}(t), with ~ \psi_{2}(t)=Q_{2T}x(T),$$
and
$$\begin{aligned}\dot{x}(t)=&Ax(t)-B_{1}R_{11}^{-1}(B_{1}+C)^{T}\psi_{1}(t)-B_{2}R_{22}^{-1}(B_{2}\\&+C)^{T}\psi_{2}(t)+CZ; x(0)=x_{0}.\end{aligned}$$
Stated differently,  if there is an open-loop Nash equilibrium, then the following differential equation has a solution:

$\frac{d}{dt}\left(
                \begin{array}{c}
                  x \\
                  \psi_{1}(t) \\
                  \psi_{2}(t) \\
                \end{array}
              \right)
=\left(
                \begin{array}{c}
                  C \\
                  0 \\
                  0 \\
                \end{array}
              \right)Z+\left(
                \begin{array}{c}
                  x \\
                  \psi_{1}(t) \\
                  \psi_{2}(t) \\
                \end{array}
              \right)$
$$\left(
   \begin{array}{ccc}
     A& -B_{1}R_{11}^{-1}(B_{1}+C)^{T} & -B_{2}R_{22}^{-1}(B_{2}+C)^{T} \\
     -Q_{1} & -A^{T} & 0 \\
     -Q_{2} & 0 & -A^{T}  \\
   \end{array}
 \right),$$
with boundary value conditions $\left\{
                                   \begin{array}{lll}
                                     x(0)=x_{0},  \\
                                     \psi_{1} (T)-Q_{1T}X(T)=0, \\
   \psi_{2}(T)-Q_{2T}X(T)=0.
                                   \end{array}
                                 \right.
$
Let $$N=\left(
                \begin{array}{c}
                  C \\
                  0 \\
                  0 \\
                \end{array}
              \right),
 y(t):=\left(
                \begin{array}{c}
                  x \\
                  \psi_{1}(t) \\
                  \psi_{2}(t) \\
                \end{array}
              \right),$$
$$
M:=\left(
   \begin{array}{ccc}
     A& -B_{1}R_{11}^{-1}(B_{1}+C)^{T} & -B_{2}R_{22}^{-1}(B_{2}+C)^{T} \\
     -Q_{1} & -A^{T} & 0 \\
     -Q_{2} & 0 & -A^{T}  \\
   \end{array}
 \right),$$
then for each $x_{0}$, the existence of a Nash equilibrium action can be transformed into the following linear system with a two-point boundary value problem having a solution:
$$\dot{y}(t)=My(t)+NZ(t),$$
with $$Py(0)+Qy(T)=\left(
                \begin{array}{c}
                  x_{0} \\
                  0 \\
                  0 \\
                \end{array}
              \right),$$
$$ P=\left(
   \begin{array}{ccc}
     I&0 &0\\
   0 & 0 & 0 \\
    0 & 0 &0  \\
   \end{array}
 \right),  Q=\left(
   \begin{array}{ccc}
     0&0 &0\\
   -Q_{1T} & I & 0 \\
    -Q_{2T} & 0 &I  \\
   \end{array}
 \right) .$$
The unique solution is
$$y(t)=e^{M(t-t_{0})}y(0)+e^{Mt}\int_{t_{0}}^{t}e^{M\tau}NZ(\tau)d\tau,$$
$$y(T)=e^{MT}y(0)+e^{MT}\int_{t_{0}}^{t}e^{M\tau}NZ(\tau)d\tau,$$
$$Py(0)+Qe^{MT}y(0)+Qe^{MT}\int_{0}^{T}e^{M\tau}NZ(\tau)d\tau=\left(
                \begin{array}{c}
                  x_{0} \\
                  0 \\
                  0 \\
                \end{array}
              \right), $$
$$(Pe^{-MT}+Q)e^{MT}y(0)=\left(
                \begin{array}{c}
                  x_{0} \\
                  0 \\
                  0 \\
                \end{array}
              \right)-Qe^{MT}\int_{0}^{T}e^{M\tau}NZ(\tau)d\tau.$$
 A necessary and sufficient condition for a unique solution to the above equation is that the matrix $(Pe^{-MT}+Q)$ is invertible, then
$$\begin{footnotesize}y(0)=e^{-MT}(Pe^{-MT}+Q)\left(
                            \begin{array}{c}
                              \left(
                                 \begin{array}{c}
                                   x_{0} \\
                                   0 \\
                                   0 \\
                                 \end{array}
                               \right)-Qe^{MT}\int_{0}^{T}e^{M\tau}NZ(\tau)d\tau
 \\
                            \end{array}
                          \right).
\end{footnotesize}$$
It follows that $y(0)$ is uniquely determined, so if there is an open-loop Nash equilibrium for all $x_{0}$, then there is a Nash equilibrium action for each $x_{0}$.
(``$\Longleftarrow$" part)
Suppose that the Riccati equation has a solution on $[0,T]$ and $(Pe^{-MT}+Q)$ is invertible. In this case, the two-point boundary value problem has a unique solution for each $x_{0}$.
Let $$y(t)= \left(
                                 \begin{array}{c}
                                   x(t) \\
                                      \psi_{1}(t) \\
                                  \psi_{2}(t)
                                 \end{array}  \right),
m_{i}(t):= \psi_{i}(t) -K_{i}(t)x(t), m_{i}(T)=0.$$

$\dot{m}_{i}(t)$
\begin{eqnarray}    \label{eq}
&=&\dot{\psi}_{i}(t)-\dot{K}_{i}(t)x(t)-K_{i}(t)\dot{x}(t)\nonumber    \\
&=&-Q_{i}x(t)-A^{T}{\psi}_{i}(t)-[-A^{T}K_{i}(t)-K_{i}(t)A \nonumber    \\
&\;&+K_{i}(t)S_{i}(t)K_{i}(t)-Q_{i}]x(t)-K_{i}(t)[Ax(t)-S_{1} \psi_{1}(t)  \nonumber    \\ &\;&-S_{2} \psi_{2}(t)+CZ(t)]\nonumber    \\  
&=& -A^{T}m_{i}(t)+K_{i}(t)Ax(t)-K_{i}(t)S_{i}K_{i}(t)x(t)
 \nonumber    \\
&\;&
-K_{i}(t)Ax(t)+K_{i}(t)S_{1} \psi_{1}(t)+K_{i}(t)S_{2} \psi_{2}(t)
\nonumber    \\
&\;&-K_{i}(t)S_{1}K_{1}(t)x(t)-K_{i}(t)S_{2}K_{2}(t)x(t)+
\nonumber    \\
&\;&
K_{i}(t)S_{1}K_{1}(t)x(t)+K_{i}(t)S_{2}K_{2}(t)x(t)-K_{i}(t)CZ(t) \nonumber    \\
&=&-A^{T}m_{i}(t)+K_{i}(t)[S_{1}m_{1}(t)+S_{2}m_{2}(t)-CZ(t)] 
\nonumber    \\
&\;&+ K_{i}(t)[-S_{i}K_{i}(t)+S_{1}K_{1}(t)+S_{2}K_{2}(t)]x(t)  \nonumber    \\
\end{eqnarray}

Now let us think about the minimum of $J_{1}(u_{1}, u_{2}^{*})$ when $u^{*}_{i}$ is equal to $-R_{ii}^{-1}(B_{i}+C)^{T}(K_{i}x+m_{i})$:
\begin{eqnarray*}
min_{u_{1}}J_{1}(u_{1}, u_{2}^{*})&=&\int_{0}^{T}\{x^{T}(t)Q_{1}x(t)+u_{1}^{T}(t)R_{11}u_{1} (t)+
\nonumber    \\
&\;&u_{2}^{*T}(t)R_{12}u_{2} (t)^{*}\}dt+x^{T}(T)Q_{1T}x(T),\end{eqnarray*}
where $\dot{x}=Ax+B_{1}u_{1}+B_{2}u_{2}^{*}+CZ, x(0)=x_{0},$  has a unique solution $\tilde{u}_{1}(t)$, $$\tilde{u}_{1}(t)=-R_{11}^{-1}(B_{1}+C)^{T}[K_{1}(t)\tilde{x}(t)+\tilde{m}_{1}(t)],$$
where $\tilde{m}_{1}(T)=0$ is the solution of a linear differential equation
$$\dot{\tilde{m}}_{1}(t)=[K_{1}(t)S_{1}-A^{T}]\tilde{m}_{1}(t)-K_{1}(t)[B_{2}u_{2}^{*}(t)+CZ(t)], $$ $$ \tilde{m}_{1}(T)=0, S_{1}=B_{1}R_{11}(B_{1}+C)^{T},$$
and $\tilde{x}(t)$ is the solution to the following differential equation
$$\dot{\tilde{x}}(t)=(A-S_{1}K_{1})\tilde{x}(t)-S_{1}\tilde{m}_{1}(t)+B_{2}u_{2}^{*}(t)+CZ(t), x(0)=x_{0}.$$

Let us substitute the value of $u_{2}^{*}(t)$ into $\dot{\tilde{m}}_{1}(t), \dot{\tilde{x}}(t)$, then
\begin{eqnarray*}
\dot{\tilde{m}}_{1}(t)&=&K_{1}(t)S_{1}\tilde{m}_{1}(t)-A^{T}\tilde{m}_{1}(t)+K_{1}(t)B_{2}R_{22}^{-1}(B_{2}
\nonumber    \\
&\;&
+C)^{T}K_{2}(t)x(t)+K_{1}(t)B_{2}R_{22}^{-1}(B_{2}+C)^{T}m_{2}\nonumber    \\
&\;&-K_{1}CZ,\end{eqnarray*}
$$\begin{aligned}\dot{m}_{1}(t)=&-A^{T}m_{1}+K_{1}S_{1}m_{1}+K_{1}S_{1}m_{2}-K_{1}S_{1}K_{1}x\\&+K_{1}S_{2}K_{2}x-K_{1}CZ,\end{aligned}$$
it follows that
\begin{eqnarray*}\dot{x}(t)&=&Ax+B_{1}\{-R_{11}^{-1}(B_{1}+C)^{T}[K_{1}(t)\tilde{x}(t)+\tilde{m}_{1}(t)]\}
\nonumber    \\
&\;&
+B_{2}\{-R_{22}^{-1}(B_{2}+C)^{T}[K_{2}(t)\tilde{x}(t)+\tilde{m}_{2}(t)]\}+CZ.\end{eqnarray*}
According to $\left\{
    \begin{array}{ll}
      \tilde{x}(t)=x(t), \\
      \tilde{m}_{1}(t)= m_{1}(t),
    \end{array}
  \right.$
 the following formula can be obtained: $$\tilde{u}_{1}(t)=-R_{11}^{-1}(B_{1}+C)^{T}[K_{1}(t)\tilde{x}(t)+m_{1}(t)].$$
And according to the uniqueness of the solution, then
$u_{1}(t)=\tilde{u}_{1}(t)=u^{*}_{1}(t).$
That is, $J_{1}(u^{*}_{1},u^{*}_{2})\leq J_{1}(u_{1},u^{*}_{2})$. Similarly, $J_{2}(u^{*}_{1},u^{*}_{2})\leq J_{2}(u^{*}_{1},u_{2})$.
Thus,  $(u_{1}^{*}, u_{2}^{*})$ is a Nash equilibrium action.
\end{proof}

When the number of players is $H$,  similarly, the Nash equilibrium action can be proved to be $u_{i}=-R_{ii}^{-1}(B_{i}+C)^{T}\varphi_{i}(t)$.

\remark Theorem 2 in \cite{zhangrenren} obtained the algebraic conditions for judging the game-based control system when regulator's control is assumed to be 0. Through the above analysis, we find that when regulator's control is not 0, Nash equilibrium action $u_{i}^{*}$  changes to $u_{i}=-R_{ii}^{-1}(B_{i}+C)^{T}\varphi_{i}(t)$. This change is reflected in the matrix $\bar{A}$ , and it can be seen from its proof that this change does not affect its algebraic conditions.

According to Theorem 1, we can obtain
\begin{equation}\label{gs34}
\left\{
  \begin{array}{ll}
    u_{i}^{*}=R_{i}^{-1}(B+C)\psi_{i}(t) \\
    \dot{ \psi}_{i}(t)=Q_{i}X(t)-\tilde{A}\psi_{i}(t) \\
     { \psi}_{i}(T)=-\tilde{Q}_{iT}X(T), i=1,2,...,H
  \end{array}
\right.
\end{equation}
Then, (\ref{gs7}) can be rewritten as follows:
\begin{equation}\label{gs35}
\left\{
  \begin{array}{ll}
     \dot{X}(t)=AX(t)+\sum_{i=1}^{H}\tilde{B}_{i}u_{i}^{*}+\tilde{B}u \\
   {X}(0)=X_{0}
  \end{array}
\right.
\end{equation}
According to (\ref{gs34}) and (\ref{gs35}), it can be obtained
\begin{equation}\label{gs36}
  \left(
    \begin{array}{c}
      \dot{X}(t) \\
      \dot{\Psi }(t)\\
    \end{array}
  \right)=\bar{A}  \left(
    \begin{array}{c}
      {X}(t) \\
      {\Psi }(t)\\
    \end{array}
  \right)+\bar{B}u.
\end{equation}

According to \cite{zhangrenren}, the controllability of system (7)  (8) is equivalent to that of system (\ref{gs36}), so the controllability of system (\ref{gs36}) is studied in the following.
In particular, the controllability matrix of the game-based control system is:
\begin{equation}\label{gs37}
 Q=\left(
     \begin{array}{ccccc}
       \bar{B} & \bar{A}\bar{B} & \bar{A}^{2}\bar{B} & \dots &  \bar{A}^{(H+1)^{2}-1}\bar{B} \\
     \end{array}
   \right).
\end{equation}

\section{Controllability analysis of graph-based GBCS}
When we use algebraic conditions to judge the controllability of game-based control systems, we need to carry out complex mathematical calculations, which prompts us to look for graph-theoretic conditions to judge the controllability, and the above algebraic conditions are the theoretical support of our graph-theoretic conditions.

According to  system (\ref{gs7}) (\ref{gs8}) , it can be seen that $u_{i}$ plays a decisive role in the system, and the Nash equilibrium action considered in this paper is composed of the coefficient $\tilde{B}_{i}$ in front of $u_{i}$. In order to study the essential influence of topology structure on game-based control system more concretely, let $\tilde{A}$ be  identity matrix, $C=0$, and combined with the particularity of $\tilde{B}_{i}$ in the system, we make the action  influence between agents reflected in whether they are neighbors.
If agent $i$ and $j$ are neighbors, it is reflected in  $\tilde{B}_{i}(j)=1,$ where $\tilde{B}_{i}(j)$ represents the $j$-th element in vector $\tilde{B}_{i}$, otherwise, $\tilde{B}_{i}(j)=0$.

\definition{Strategy equivalence partition (SEP):  If different agents receive the same number of strategies from agents in any cell, then these agents can be divided into the same cell.
Let $s(i)$ denote the number of strategies of agent $i$, and $C_{p}$ denote the $p$-th cell,  then for any $i, j, p, (s(i), C_{p})=(s(j), C_{p})$},  agents $i,j$ can be divided into the same cell, where $(s(j), C_{p})$   represents the total number of strategies of the cell $C_{p}$ received by agent  $j$. If the number of nodes in a cell is greater than one, it is called nontrivial, otherwise it is called trivial.

\theorem{Suppose that Assumption 1 hold.  If  there is a nontrivial strategy equivalent partition, $C_{i}=\{i_{1},i_{2},...,i_{p}\}$, and  the $i_{1}$-th ,$i_{2}$-th ,...,$i_{q}$-th row vectors of $T$ are equal, where
$$ T=[I_{(H+1)}~0]e^{-\bar{A}T}\left(
                                     \begin{array}{c}
                                       I_{H+1} \\
                                       -\tilde{Q}_{1T} \\
                                       \vdots \\
                                        -\tilde{Q}_{HT} \\
                                     \end{array}
                                   \right)
\left(
                                                     \begin{array}{c}
                                                       0_{1\times(H)} \\
                                                       I_{H} \\
                                                     \end{array}
                                                   \right),$$
then, the GBCS (\ref{gs36}) is uncontrollable.
Moreover, the general formula of controllability matrix (\ref{gs37})   is
$Q_{pq} =$

$$\left\{
           \begin{array}{ll}
             Q_{1q}=\left\{
                           \begin{array}{ll}
                             \tilde{B}, & \hbox{$q=1$;} \\
                           \begin{aligned} &\tilde{A}^{2}Q_{1(q-2)}+
\\&\sum_{i=1}^{H} \tilde{B}_{i}R_{i}^{-1} \tilde{B}_{i}^{T}Q_{i}Q_{1(q-2)} \end{aligned}, & \hbox{$q$ is  odd;} \\
                              \tilde{A}Q_{1(q-1)}, & \hbox{$q$ is even.}
                           \end{array}
                         \right.
\\
           \begin{aligned}   & Q_{pq}\\& (1<p \leq (H+1)^{2}-1) \end{aligned}=\left\{
                                            \begin{array}{ll}
                                              0, & \hbox{$q$ is  odd;} \\
                                              Q_{(p-1)} Q_{1(q-1)}, & \hbox{$q$ is even.}
                                            \end{array}
                                          \right.
           \end{array}
         \right.$$}

\begin{proof}
To prove the controllability condition of the game-based control system  (\ref{gs36}), we first prove the general formula of controllability matrix. We first compute the first element of the three cases in $Q_{1q}$:
$$Q_{11}=\tilde{B}, Q_{12}=\tilde{A}\tilde{B}, Q_{13}=\tilde{A}^{2}\tilde{B}+\sum_{i=1}^{H}\tilde{B}_{i}R_{i}^{-1}(\tilde{B}_{i}+C)^{T}Q_{i}\tilde{B}.$$

(1) If $(H+1)^{2}$ is even, then $(H+1)^{2}-1$ is odd, that is, the last element in $Q_{1q}$ is odd, then we assume that
$$  \begin{aligned} Q_{1[(H+1)^2-3]}=&  \tilde{A}^{2}Q_{1[(H+1)^2-5]}\\&+\sum_{i=1}^{H} \tilde{B}_{i}R_{i}^{-1} (\tilde{B}_{i}+C)^{T}Q_{i}Q_{1[(H+1)^2-5]},  \end{aligned}$$ $$Q_{1[(H+1)^2-4]}=\tilde{A}Q_{1[(H+1)^2-5]},$$
the calculations show that
$$Q_{21}=...=Q_{H1}=0, Q_{i2}=\tilde{A}Q_{i}, i=2,...,H,$$ then
$$Q_{23}=Q_{1}Q_{12}-\tilde{A}^{T}Q_{1}\tilde{B}.$$
And since $Q_{1}$ is a symmetric matrix,  $Q_{23}=0$. Similarly, $Q_{33}=Q_{43}=...=Q_{H3}=0$. Suppose that $$Q_{2[(H+1)^{2}-3]}=Q_{3[(H+1)^{2}-3]}=...=Q_{H[(H+1)^{2}-3]}=0,$$
$$Q_{i[(H+1)^{2}-4]}=Q_{(i-1)}Q_{1[(H+1)^{2}-5]}, i=2,3,...,H.$$
then, $$
Q_{1[(H+1)^{2}-2]}=\tilde{A}Q_{1[(H+1)^{2}-3]}+\sum_{i=1}^{H} \tilde{B}_{i}R_{i}^{-1} (\tilde{B}_{i}+C)^{T}\times0,
$$
$$Q_{i[(H+1)^{2}-2]}=Q_{i}Q_{1[(H+1)^{2}-3]}+(-\tilde{A}^{T})\times0,$$
$$\scriptsize Q_{1[(H+1)^{2}-1]}=\tilde{A}Q_{1[(H+1)^{2}-2]}+\sum_{i=1}^{H} \tilde{B}_{i}R_{i}^{-1} (\tilde{B}_{i}+C)^{T}Q_{i}Q_{1[(H+1)^{2}-3]},$$
$$Q_{i[(H+1)^{2}-1]}=Q_{i}\tilde{A}Q_{1[(H+1)^{2}-3]}-\tilde{A}^{T}Q_{i}Q_{1[(H+1)^{2}-3]},$$

(2) If $(H+1)^{2}$ is odd, the last column of the controllability matrix   is even, then assume that $$Q_{1[(H+1)^2-3]}= \tilde{A}Q_{1[(H+1)^2-4]} , $$ $$Q_{1[(H+1)^2-4]}=\tilde{A}^{2}Q_{1[(H+1)^2-6]}+\sum_{i=1}^{H} \tilde{B}_{i}R_{i}^{-1} (\tilde{B}_{i}+C)^{T}Q_{i}Q_{1[(H+1)^2-6]}.
$$
Its initial elements do not change. In this case, the difference from (1) is the subscript of the element, and the proof details  are omitted.

According to the above general formula of controllability matrix, it can be found that matrix  $\sum_{i=1}^{H} \tilde{B}_{i}R_{i}^{-1} \tilde{B}_{i}^{T}$  exists in every element of controllability matrix. To obtain a clearer view of this question within topology structure, take
$R_{i}^{-1}=1$. The following focuses on the properties of matrix $\sum_{i=1}^{H} \tilde{B}_{i}R_{i}^{-1} \tilde{B}_{i}^{T}$ in the presence of strategy equivalent partition.
 Let $S=\sum_{i=1}^{H} \tilde{B}_{i} \tilde{B}_{i}^{T}, S\in R^{(H+1)\times (H+1)}$,
 set the serial number of regulator in the system as $r_{1}$, and the serial numbers of other agents as $r_{2},r_{3},...,r_{H}.$
Regulator has a neighbor relationship with each agent, and it can also be called neighbor relationship with itself, so the non-neighbor row of vector $\tilde{B}_{i}$ is 0. That is, the nonzero elements in matrix $\tilde{B}_{i} \tilde{B}_{i}^{T}$ are located in  the neighbor rows and neighbor columns of agent $i$.
 Therefore, the element in row $i$ and column $j$  of matrix $\sum_{i=1}^{H} \tilde{B}_{i} \tilde{B}_{i}^{T}$ represents the number of agents containing both $i$ and $j$ in the neighbor set, that is, $S_{ij}$ is equal to the number of common neighbors of  $i$ and $j$.
 It is worth noting that $S_{ii}$ denotes the number of neighbors of agent $i$, including $i$ itself. $S$  must be a symmetric matrix, and the first row and column are equal to the diagonal elements.
 If there is a strategy equivalent partition in the system, it is assumed that there are $p$ cells, $C_{1}, ... ,C_{p}.$  If $|C_{i}|=s_{i}$, that is,  the cell $C_{i}$ contains $s_{i}$  elements. The nodes in the cell  are marked with $s_{i}$  consecutive numbers, and then the strategy matrix $S$ is divided into blocks according to the order of the cell:
$$S=\left(
    \begin{array}{cccc}
      P_{11} &  P_{12} & \dots &  P_{1p} \\
       P_{21} & P_{22} &   \dots& P_{2p} \\
       \vdots & \vdots &   &  \vdots \\
       P_{p1} & P_{p2} & \dots & P_{pp} \\
    \end{array}
  \right),
$$
where $P_{ij} \in R^{s_{i}\times s_{j}}$.  Let $i_{1}, , i_{2} \dots , i_{q} \in C_{i}$,  then
$$\sum_{j=1}^{(H+1)}S_{i_{1}j}=\sum_{j=1}^{(H+1)}S_{i_{2}j}=\dots =\sum_{j=1}^{(H+1)}S_{i_{n}j}.$$
That means the sum of the elements in the same cell is equal. Or, stated differently, lines $i_{1},\dots, i_{q}$  in   $\sum_{i=1}^{H} \tilde{B}_{i}R_{i}^{-1} \tilde{B}_{i}^{T} \tilde{B}$  are equal.
Therefore, in the matrix
$$[I_{(H+1)}~0]\bar{A}^{k+1}\bar{B}, k=1,2,...,[(H+1)^2-1].$$
the rows of the above matrix corresponding to the elements of the same cell are equal. If the entries in the same cell correspond to the row vectors of the following matrix  are also equal:
$$\scriptsize[I_{(H+1)}~0]e^{-\bar{A}T}\left(
                                     \begin{array}{c}
                                       I_{H+1} \\
                                       -\tilde{Q}_{1T} \\
                                       \vdots \\
                                        -\tilde{Q}_{HT} \\
                                     \end{array}
                                   \right)
\left(
                                                     \begin{array}{c}
                                                       0_{1\times(H)} \\
                                                       I_{H} \\
                                                     \end{array}
                                                   \right),$$
then the  following matrix  is not full row rank:
$$\scriptsize \left(
    \begin{array}{cc}
       [I_{(H+1)}~0]\bar{A}^{k+1}\bar{B} & [I_{(H+1)}~0]e^{-\bar{A}T}\left(
                                     \begin{array}{c}
                                       I_{H+1} \\
                                       -\tilde{Q}_{1T} \\
                                       \vdots \\
                                        -\tilde{Q}_{HT} \\
                                     \end{array}
                                   \right)
\left(
                                                     \begin{array}{c}
                                                       0_{1\times(H)} \\
                                                       I_{H} \\
                                                     \end{array}
                                                   \right) \\
    \end{array}
  \right).
$$
According to (\cite{zhangrenren}Theorem 2), the game-based control system is  uncontrollable.
\end{proof}

\begin{figure}[!h]
  \centering
  \includegraphics[width=1.3in]{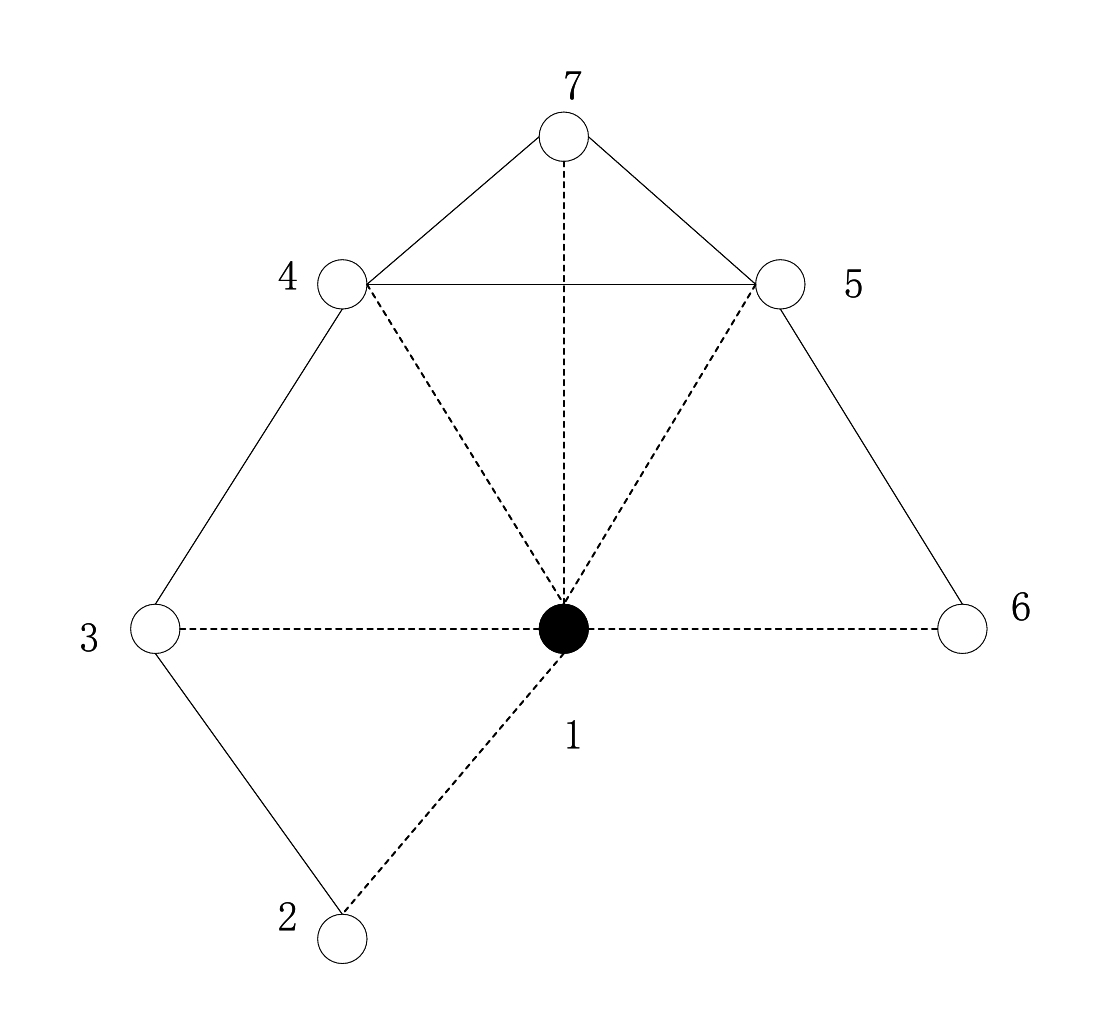}\\
  \caption{A system with 1 regulator and 6 agents}\label{F1}
\end{figure}
In control systems that do not consider game factors,   \cite{wo} shows that the equivalence partition has limitations, that is to say, whether the system is controllable can not be judged solely based on the equivalent partition. Under the model of \cite{wo}, the existence of equivalent partition can definitely lead to  uncontrollable system, but  uncontrollable sysem is not necessarily because of the existence of equivalent partition.
As shown in Figure 1, node 1 is considered as the leader and the remaining nodes are selected as followers.
It is judged that there is no equivalent partition in the system, but the system  is still  uncontrollable, because the system matrix has eigenvectors orthogonal to $1_{n}$ vector. However, under the game-based control system, the system composed of this topology structure is controllable. This is the only type of counterexample due to the limitation of the equivalence partition problem,  but this structure does not hold for game-based control systems, because the strategy matrix $S$ does not have eigenvectors that are orthogonal to the $1_{n}$ vector. Therefore, we make a conjecture that there is no limitation of strategy equivalence partition in game-based control systems, that is, the sufficiency of Theorem 2 holds.

\section{Conclusion}
What effect does topology have on the game-based control system? This paper holds that its essential influence lies in the Nash equilibrium action, that is to say, different topologies will produce different Nash equilibrium actions. We first considered the expression of Nash equilibrium action under the premise that regulator's control strategy is not zero.  Then, we obtained the general formula of controllability matrix of game-based control system. Based on this, we found the "Laplacian matrix in game-based control system "--the strategy matrix, which can not only be obtained by algebraic calculation, but also can be directly written by topological structure. Perhaps the position of strategy matrix in game-based control system is comparable to that of Laplacian matrix in general control system. Through the above analysis,  we finally obtained the graph theory condition based on the strategy equivalent partition, and put forward the conjecture that the strategy equivalent partition does not have the limitation in the GBCS.

\end{document}